\documentclass[reqno]{amsart}
\usepackage{amscd} 
\usepackage{amsfonts} 
\usepackage{amssymb} 
\usepackage{latexsym} 
\input{diagrams}
\newcommand{\ncm}{\newcommand}

\def\ME{\mathcal{E}}

\newtheorem{theorem}{Theorem}[section]
\newtheorem{prop}[theorem]{Proposition}
\newtheorem{lemma}[theorem]{Lemma}
\newtheorem{cor}[theorem]{Corollary}

\newtheorem{lem&def}[theorem]{Lemma \& Definition}

\newtheorem{remark}[theorem]{Remark}

\ncm{\End}{\mbox{\rm End}\,}

\def\Hom{\mbox{\rm Hom}\,}

\def\id{\mbox{\rm id}}

\def\into{\hookrightarrow}
\def\to{\rightarrow}

\def\o{\otimes}    

\ncm{\rarr}[1]{\stackrel{#1}{\longrightarrow}}
\ncm{\larr}[1]{\stackrel{#1}{\longleftarrow}}
\def\cop{\Delta}

\def\eps{\varepsilon}

\def\du1{\hat 1}

\def\-1{_{(-1)}}
\def\0{_{(0)}}
\def\1{_{(1)}}
\def\2{_{(2)}}
\def\3{_{(3)}}
\def\4{_{(4)}}
\def\I{^{(I)}}

\def\|{\, | \,}

\def\du1{\hat 1}

\def\ract{\triangleleft}

\begin{document}

\title[Depth two and infinite index]{Infinite index subalgebras of depth
two}
\author{Lars Kadison}
\address{Department of Mathematics \\ 
University of Pennsylvania \\ 209 South 33rd Street \\ Philadelphia, PA 19104-6395} 
\email{lkadison@c2i.net} 
\date{}
\thanks{}
\subjclass{16W30 (46L37, 81R50)}  

\begin{abstract} 
An algebra extension $A \| B$ is right depth two in this paper
if its tensor-square is $A$-$B$-isomorphic to a direct summand
of any (not necessarily finite) direct sum of $A$ with itself.
For example, normal subgroups of infinite groups, infinitely generated
Hopf-Galois extensions and infinite dimensional
algebras are
 depth two in this extended sense.
The added generality loses some duality results obtained in the finite theory \cite{KS} but  extends the main theorem of depth two theory, as for example in \cite{KN}. 
That is, a right depth two extension has right bialgebroid $T = (A \o_B A)^B$ over its centralizer
$R = C_A(B)$. The main theorem:  an extension $A \| B$
is right depth two and right balanced if
and only if $A \| B$
is $T$-Galois wrt.\ left projective, right $R$-bialgebroid $T$.         
\end{abstract} 
\maketitle

\section{Introduction}

Bialgebroids arise as the endomorphisms of fiber functors from certain tensor categories to a bimodules category over a base algebra.  For example, bialgebras are bialgebroids over a one-dimensional base algebra, while weak bialgebras are bialgebroids over a separable base algebra. Hopf algebroids are  bialgebroids with antipodes: various twisted Hopf algebras are also Hopf algebroids over a one-dimensional base algebra.  

Like bialgebras and their actions/coactions, bialgebroids also act and coact on noncommutative algebras in a more general setting suitable to mathematical physics \cite{BB}. Initially
appearing in the analytic theory of subfactors, the notion of depth two has been widened to Frobenius extensions in
\cite{KN} and to arbitrary subalgebras
in \cite{KS}. As shown in \cite{KS}
and later papers, depth two is a \textit{Galois} theory of actions and coactions for 
bialgebroids.  In this paper, we widen the definition of depth
two algebra extension in \cite{KS} to include Hopf $H$-Galois
extensions where $ H$ is an infinite-dimensional Hopf algebra, 
such as the univ.\ env.\ algebra of a Lie
algebra or an infinite dimensional group
algebra.  Although
we lose the dual theory of finite projective, left and right bialgebroids
over the centralizer in \cite{KS}, we retain the right
bialgebroid $T$ and its role in coaction in \cite{LK2005}.  We  then obtain the main theorem of depth
two Galois theory with no finiteness conditions (Theorem~\ref{th-main}):  an algebra
extension $A \| B$ is right depth two with $A_B$ 
a balanced module if and only if
$A \| B$ is $T$-Galois w.r.t.\
a left projective right $R$-bialgebroid $T$, for some base ring $R$
which  commutes within $A$
with the subring of 
coinvariants $B$.

\subsection{Depth two preliminaries}
By algebra we mean a unital associative
algebra over a commutative ring $k$,
and by algebra extension $A \| B$, we mean
any identity-preserving algebra homomorphism
$B \to A$, proper if $B \to A$ is monic.
In either case,   the natural
bimodule ${}_BA_B$ and its properties
define the properties of the extension
from this point of view. For example,
we say $A \| B$ is right faithfully flat 
if $A_B$ is faithfully flat, in which
case one notes the
extension $A \| B$ is proper.  

An algebra extension $A \| B$ is \textit{left depth two (D2)} if
its tensor-square $A \o_B A$ as a natural $B$-$A$-bimodule
is isomorphic to a direct summand of a direct sum
of the natural $B$-$A$-bimodule $A$: equivalently, for some set 
$I$, we have 
\begin{equation}
\label{eq: D2}
A \o_B A \oplus * \cong A^{(I)},
\end{equation}
where $A^{(I)}$ denotes the coproduct (weak direct
product, direct sum $\sum_{i \in I} A_i$,
each $A_i = A$) of $A$ with itself indexed by
$I$ and consists of elements $(a_i)_{i\in I}$
where $a_i \in A$ and $a_i = 0$ for all but finitely many indices
(almost everywhere, a.e.).  
An extension $A \| B$ is \textit{right D2} if eq.~(\ref{eq: D2}) holds instead as natural
$A$-$B$-bimodules. An algebra extension is of course D2 if it is both left D2 and right D2. 

For example, 
if $A \| B$ is a projective algebra
(so $B$ is commutative, maps into the center of $A$
and the module $A_B$ is projective), then
$A \| B$ is D2, since $A_B \oplus * \cong
B^{(I)}$ for index set $I$, so
we may tensor this by $- \o_B A_A$
to obtain eq.~(\ref{eq: D2}).   

As another example, suppose $H$ is a Hopf algebra
of finite or infinite dimension over a field,
and $A$ is a right $H$-comodule algebra with
$B$ equal to the subalgebra of coinvariants.
If $A \| B$ is an $H$-Galois extension,
then $A \| B$ is right D2 since $A \o_B A \cong
A \o H$ via the Galois $A$-$B$-isomorphism,
$x \o y \mapsto x y\0 \o y\1$, where $y\0 \o y\1$
denote finite sums of elements equal to the 
value in $A \o H$ of the coaction on $y \in A$.  Let $I$ be
in one-to-one correspondence with a basis for $H$.
Then $A \o_B A \cong A^{(I)}$. If $H$ has a bijective
antipode,  use the equivalent Galois $B$-$A$-bimodule isomorphism
given by $x \o y \mapsto x\0 y \o x\1$ to conclude that $A \| B$
is left D2.    

If the index set $I$ is finite, then
the algebra extension 
$A \| B$ is right or left D2 in the earlier 
sense of \cite{KN, KS, LK2003, LK2005}.  
The lemma below notes that  the
earlier definition is recovered for any f.g.\
extension.  
\begin{lemma}
If $A \| B$ is right  or left D2 and either of the natural modules ${}_BA$ or $A_B$ is
finitely generated, then $I$
in eq.~(\ref{eq: D2}) may be
chosen finite.  
\end{lemma}
\begin{proof}
Suppose $A \| B$ is right D2.  If either
${}_BA$ or $A_B$ is f.g., then ${}_AA\o_B A_B$ is f.g.
It follows that ${}_AA \o_B A_B$
is isomorphic to a direct summand
of a finite direct sum $A^n \subseteq A^{(I)}$.
The argument is entirely similar starting
with a left D2, left or right f.g.\ extension.

More explicitly using the $A$-$B$-epi $f$
and $A$-$B$-monic $g$ defined below, if $A \o_B A = 
Aw_1 B + \cdots + Aw_N B$ for $N$ elements
$w_j \in A \o_B A$, then 
$g(w_j) = (a_{ij})_{i \in I}$
has finite support on $I_j \subset I$,
then $I' = I_1 \cup \cdots \cup I_N$
is finite and $g$ corestricts, $f$
restricts to $A^{I'}$ so that $f \circ g
= \id_{A \o_BA}$.  
\end{proof}

In analogy with projective bases for projective modules,
we similarly develop D2 quasibases for depth two extensions.

\begin{prop}
An algebra extension is right D2 if and only if
there is an index set $I$ and sets of elements $u_i
= u^1_i \o_B u^2_i \in (A \o_B A)^B$, $\gamma_i \in \End {}_BA_B$, both indexed by $I$, such that for each
$a \in A$, $\gamma_i(a) = 0$ a.e.\ on $I$,  and 
\begin{equation}
\label{eq: rd2qb}
x \o_B y = \sum_{i \in I} x\gamma_i(y)u^1_i \o_B u^2_i
\end{equation}
for all
$x,y \in A$.
\end{prop}
\begin{proof}
Let $\pi_i: A^{(I)} \to A$ and $\iota_i: A \to A^{(I)}$ be the usual projection and inclusion mappings of a coproduct, so that $\pi_j \circ \iota_i
= \delta_{ij} \id_A$ and $\sum_{i \in I} \iota_i \circ \pi_i = \id$ on $A^{(I)}$. 

Given a right D2 extension $A \| B$, there is
an $A$-$B$-split epimorphism $f: A\I \to A\o_B A$,
say with section $g : A\o_B A \to A\I$.
Then $f \circ g = \id_{A \o_B A}$.
Define $f_i = f \circ \iota_i \in \Hom (A, A \o_B A)$
and define $g_i = \pi_i \circ g \in \Hom (A \o_B A, A)$, both  hom-groups of the natural $A$-$B$-bimodules. Then
$\sum_{i \in I} f_i \circ g_i = \id_{A \o_B A}$.   But
\begin{equation}
\Hom (A, A \o_B A) \cong (A \o_B A)^B
\end{equation}
via $f \mapsto f(1_A)$, and 
\begin{equation}
\Hom (A \o_B A, A) \cong \End {}_BA_B
\end{equation}
via $F \mapsto F(1_A \o -)$ with inverse
$\alpha \mapsto (x \o y \mapsto x\alpha(y))$.  
In this case, there are $\gamma_i \in \End {}_BA_B$
such that $\gamma_i(a) = g_i(1 \o a)$, all $a \in A$, and $u_i \in (A \o_B A)^B$ such that
$f_i(1_A) = u_i$, for each $i \in I$.
Note that $\gamma_i(a) = 0$ a.e.\ on $I$,
since $g_i(1 \o a) = \pi_i(g(1 \o a))$ is
zero a.e.\ on $I$.    
It follows from $\id_{A \o_BA} = \sum_{i \in I} f_i \circ g_i$ that
$$ x \o y = \sum_{i \in I} x\gamma_i(y)u_i. $$

Conversely, given right D2 quasibases $\{ \gamma_i \}_{i \in I}$, $\{ u_i \}_{i \in I}$ as above,
define epimorphism $\pi: A\I \to A \o_B A$ 
of natural $A$-$B$-bimodules by
\begin{equation}
\pi[ (a_i) ] = \sum_{i \in I} a_i u_i
\end{equation}
with $A$-$B$-bimodule section $\iota: A\o_B A \into A\I$ given by  
\begin{equation}
\iota(x \o y) = (x\gamma_i(y))_{i \in I}
\end{equation}
well-defined in $A\I$ since for 
all $a \in A$, $\gamma_i(a) = 0$ a.e.\ on $I$.  
\end{proof}

A similar proposition holds for a left D2 extension $A \| B$
and left D2 quasibase
$t_i \in (A \o_B A)^B$ and $\beta_i \in \End {}_BA_B$
for each $i \in I$.  In this case,
\begin{equation}
\label{eq: ld2qb}
\sum_{i \in I} t^1_i \o_B t^2_i\beta_i(x)y =   x \o_B y, 
\end{equation}
for all $x,y \in A$, 
  which is equivalently expressed as $a \o 1 = \sum_i t_i \beta_i(a)$
for all $a \in A$,
where again $\beta_i(a) = 0$ a.e.\ on the index set $I$.  
We fix our notation for right and left D2 quasibases
throughout the paper.  In addition, we
denote $T = (A \o_B A)^B$ and (less importantly) $S = \End {}_BA_B$.  

For example, left and right D2 quasibases are obtained as follows
for  group algebras $A = k[G] \supseteq B = k[N]$ where $G$ is a group, possibly of infinite order,
$N$ is a normal subgroup of possibly infinite index, and 
$k$ is a commutative ring.
Let $\{ g_i \}_{i \in I}$ be a transversal of $N$ in $G$.
Define straightforwardly a projection onto the $i$'th coset by $\gamma_i(a) = \sum_{j \in J} \lambda_{ij} g_i n_j$ where $a \in A$ and therefore of the form
$a = \sum_{i \in I} \sum_{j \in J} \lambda_{ij} g_i n_j$,
where $J$ is an indexing set in one-to-one correspondence with $N$
and $k$-coefficients $\lambda_{ij} = 0$ a.e. on $I \times J$. 
In this case for any basis element $g \in G \into A$ all but
one of the projections $\gamma_i$ vanish on $g$: if $g$ is in the coset
$Ng_j$, then $\gamma_j(g) = g$.   Of course, the $\gamma_i$
are $B$-$B$-bimodule projections since $gN = Ng$ for all $g \in G$.  
It is then easy to see that
\begin{equation}
1 \o_B g = \sum_{i \in I} \gamma_i(g) g_i^{-1} \o_B g_i
\end{equation}
whence eq.~(\ref{eq: rd2qb}) follows by choosing $u_i = g_i^{-1} \o_B g_i$. Note that $u_i \in (A \o_B A)^B$
since $ng_i^{-1} \o_B g_i = g_i^{-1} \o_B g_i n$
for $n \in N$.  

Similarly a left D2 quasibase is given by $\{ \gamma_i \}$
and $\{ g_i \o_B g_i^{-1} \}$ since $g_i N = N g_i$.  

We end this section with a proposition
collecting various necessary conditions
on a right depth two algebra extension.

\begin{prop}
Suppose an algebra extension
$A \| B$ is right D2 with centralizer $R$.  Then the following is true: 
\begin{enumerate}
\item  for each two-sided ideal in $A$, $(I \cap R)A \subseteq A(I \cap R)$;
\item  ${}_BA$ is projective
if $A \| B$ is moreover a split extension.
\item for some indexing set $I$,
$\End {}_BA \oplus * \cong A^I$
as natural $B$-$A$-bimodules. 
\item  For each H-separable extension $B \| C$, or equivalently an extension satisfying
\begin{equation}
\label{eq: extH-sep}
B \o_C B \oplus \, *  \cong  B^{(J)} \ 
\mbox{\rm natural $B$-$B$-bimodules}
\end{equation}
for any index set $J$, the
composite algebra extension $A \| C$
is right D2. 
\end{enumerate}
\end{prop}
\begin{proof}
The proof of each statement follows in the order above.
\begin{enumerate}
\item Given $x \in I \cap R$ and $a \in A$, apply
 eq.~(\ref{eq: rd2qb}) and a right D2 quasibase:  $xa = \sum_i \gamma_i(a) u^1_i x u^2_i$.  
Note that for each $i$, $u^1_i x u^2_i
\in I \cap R$.  
\item Given a $B$-$B$-bimodule projection
$p: A \to B$, apply $p \o_B \id_A$ to
eq.~(\ref{eq: rd2qb}) with $x = 1$, obtaining
$y = \sum_{i \in I} p(\gamma_i(y)u^1_i)u^2_i$ for all $y \in A$, which shows
${}_BA$ has dual bases.
\item Note that $\Hom ({}_AA \o_B A, {}_AA) \cong \End {}_BA$ as $A$-$A$-bimodules
via $F \mapsto F(1 \o -)$. Apply
$\Hom ({}_A - ,{}_A A)$ to
${}_AA \o_B A_B \oplus * \cong {}_AA^{(I)}_B$, noting that
$\Hom (A^{(I)}, A) \cong A^I$ (the
direct product)
as $B$-$A$-bimodules. 

\item Apply the functor $A \o_B - \o_B A$ from
$B$-$B$-bimodules into $A$-$B$-bimodules
to the isomorphism~(\ref{eq: extH-sep}).
Then ${}_AA \o_C A_B \oplus \, * \, \cong
{}_AA \o_B A_B^{(J)}$. Clearly 
$A \o_B A^{(J)} \oplus * \cong (A^{(I)})^{(J)} \cong A^{(I \times J)}$
as $A$-$B$-bimodules. Whence
 the composite extension
$A \| C$ satisfies the right D2 condition  \begin{equation}
\label{eq: same}
A \o_C A \oplus \, * \, \cong A^{(I \times J)}.
\end{equation}
Finally, $J$ in eq.~(\ref{eq: extH-sep}) may be replaced by the finite support of the image of $1 \otimes 1$ in $A^{(J)}$, under a split $A$-$A$-monomorphism $A \o_B A \to A^{(J)}$. Whence an algebra extension
satisfying eq.~(\ref{eq: extH-sep}) is
 H-separable \cite{LK2003}. 
\end{enumerate}
\end{proof}
 Similar statements hold for a left D2 extension, one of which results in
\begin{cor}
If $A \| B$ is D2, then the centralizer
$R$ is a normal subalgebra: i.e.,
for each two-sided ideal $I$ in $A$,
the contraction of $I$ to $R$ is $A$-invariant: 
\begin{equation}
A(I \cap R) = (I \cap R)A
\end{equation}
\end{cor}

For example, any trivial extension $A \| A$ is D2, in which case
$R$ is the center of $A$, which is of course a normal subalgebra.  

\section{The bialgebroid $T$ for a depth two extension}

In this section we establish that if $A \| B$ is a right or left
D2 algebra extension, then the construct $T = (A \o_B A)^B$,
whose acquaintance we made in the last section, is a right bialgebroid
over the centralizer $C_A(B) = R$.  Moreover, $T$ is right or left
projective as a module over $R$ according to which depth two condition,
left or right, respectively, we assume.

\begin{lemma}
Let $T$ be equipped with the  natural $R$-$R$-bimodule structure given by 
\begin{equation}
\label{eq: Rbimod}
r \cdot t \cdot s = rt^1 \o_B t^2 s
\end{equation}
for each $r,s \in R$ and $t \in T$.
If $A \| B$ is left D2 (right D2), then $T$ is a projective
right (left, resp.) $R$-module.
\end{lemma}
\begin{proof}
This follows from eq.~(\ref{eq: ld2qb}) by restricting to
elements of $T \subseteq A \o_B A$.  We obtain
$t = \sum_i t_i \beta_i(t^1)t^2$ where $t_i \in T$. But $\beta_i(t^1) t^2 \in R$
so define elements $f_i \in \Hom (T_R, R_R)$, indexed by $I$, by
$f_i(t) = \beta_i(t^1)t^2$. Substitution yields $t = \sum_i t_i f_i(t)$,
where $f_i(t) = 0$ a.e.\ on $I$. Whence $T_R$ is projective
with dual basis $\{ t_i \}$, $\{ f_i \}$.  

The proof that $A \| B$ is right D2 implies ${}_RT$ is projective follows similarly from eq.~(\ref{eq: rd2qb}).
\end{proof}

The next theorem may be viewed as a generalization of the first statement in \cite[theorem 5.2]{KS}.

\begin{theorem}
\label{th-bi}
If $A \| B$ is right D2 or left D2, then $T$ is a right
bialgebroid over the centralizer $R$.  
\end{theorem}
\begin{proof}
The algebra structure on $T$ comes from the isomorphism
 $T \cong \End {}_AA \o_B A_A$ via $$t \longmapsto
(x \o_B y \mapsto xt^1 \o_B t^2 y)$$
with inverse $F \mapsto F(1 \o 1)$. The endomorphism algebra structure
on $T$ becomes 
\begin{equation}
\label{eq: tee}
tu = u^1 t^1 \o_B t^2 u^2,
\ \ \  1_T = 1_A \o_B 1_A.
\end{equation}

It follows from this that there is algebra homomorphism $s_R: R \to T$
and algebra anti-homomorphism $t_R: R \to T$, satisfying a commutativity
condition and inducing an $R$-$R$-bimodule from the right of $T$, given
by ($r,s \in R, t \in T$)
\begin{eqnarray}
s_R(r) & = & 1_A \o_B r \\
t_R(s) & = & s \o_B 1_A \\
s_R (r) t_R(s) & = & t_R(s) s_R(r) = r \o_B s \\
t t_R(r) s_R(s) & = &  rt^1 \o t^2 s. 
\end{eqnarray}
Henceforth, the bimodule ${}_RT_R$  refered to is 
the one above, which is the same as the bimodule in eq.~(\ref{eq: Rbimod}).  

An $R$-coring structure $(T,R,\cop, \eps)$ with comultiplication
$\cop: T \to T \o_R T$ and counit $\eps: T \to R$ is
given by
\begin{eqnarray}
\cop(t) & = & \sum_{i \in I} (t^1 \o_B \gamma_i(t^2)) \o_R u_i \\
\eps(t) & = & \sum t^1 t^2
\end{eqnarray}
i.e., $\eps$ is the restriction of $\mu: A \o_B A \to A$ 
defined by $\mu(x \o y) = xy$ to $T \subseteq A \o_B A$. 
The coproduct $\cop$ is well-defined since for any given
$t \in T$, there are only finitely many nonzero terms on the right.   
It is immediate that $\cop$ is left $R$-linear, $\eps$ is left
and right $R$-linear, and $$(\eps \o_R \id_T) \circ \cop = 
\id_T = (\id_T \o_R \eps) \circ \cop $$
follows from variants of eq.~(\ref{eq: rd2qb}). 
We postpone the proof of coassociativity of $\cop$
for one paragraph. 

Additionally, note that the coproduct and counit are
unit-preserving, $\eps(1_T) = 1_A = 1_R$
and $\cop(1_T) = 1_T \o_R 1_T$, since $\gamma_i(1_A) \in R$.   

We employ the usual Sweedler notation $\cop(t) = t\1 \o_R t\2$.  
In order to show the bialgebroid identities 
\begin{eqnarray}
\label{eq: rightRlin}
\cop(tr) & = & t\1 \o_R t\2 r \\
\label{eq: timesR}
s_R(r) t\1 \o_R t\2 & = & t\1 \o_R t_R(r)t\2 \\
\label{eq: homo}
\cop(tu) & = & t\1 u\1 \o_R t\2 u\2 
\end{eqnarray}
it will be useful to know that 
\begin{equation}
T \o_R T \stackrel{\cong}{\longrightarrow} (A \o_B A \o_B A)^B \ \ \ 
t \o_R u \mapsto t^1 \o t^2u^1 \o u^2
\end{equation}
with inverse $$v \mapsto \sum_i (v^1 \o_B v^2 \gamma_i(v^3)) \o_R u_i. $$

Note that the LHS and RHS of eq.~(\ref{eq: rightRlin}) are
the expressions $\sum_i (t^1 \o \gamma_i(t^2r)) \o u_i$
and $\sum_i (t^1 \o \gamma_i(t^2)) \o u_ir$, which both map bijectively into
$t^1 \o 1_A \o t^2 r$ in $(A \o_B A \o_B A)^B$, whence LHS = RHS
indeed.
Similarly, the LHS of eq.~(\ref{eq: timesR}) is
$\sum_i (t^1 \o r\gamma_i(t^2)) \o_R u_i$
while the RHS is $\sum_i (t^1 \o \gamma_i(t^2)) \o_R (u^1_i r \o u^2_i)$,
both mapping  into the same element, $t^1 \o_B r \o_B t^2 $.  
Hence, this equation holds, giving meaning to the next equation for all $t,u \in T$ (the tensor  product algebra over noncommutative
rings ordinarily makes no sense, cf.\ \cite{BW}).
The eq.~(\ref{eq: homo}) holds because both expressions map isomorphically into
the element $u^1 t^1 \o_B 1_A \o_B t^2 u^2$.  

Finally the coproduct is coassociative, $(\cop \o_R \id_T) \circ \cop = 
(\id_T \o_R \cop) \circ \cop   $ since we first note that $$T \o_R T \o_R T \stackrel{\cong}{\longrightarrow}
(A \o_B A \o_B A \o_B A)^B$$ $$t \o u \o v \mapsto
t^1 \o t^2 u^1 \o u^2 v^1 \o v^2.$$
Secondly, $\sum_i \cop(t^1 \o \gamma_i(t^2)) \o_R u_i$
maps into $t^1 \o 1 \o 1 \o t^2$,
as does $\sum_i (t^1 \o \gamma_i(t^2)) \o_R \cop(u_i)$,
which establishes this, the last of the axioms
of a right bialgebroid.

The proof that $T$ is a right bialgebroid using a left D2 quasibase instead
is very similar.  
\end{proof}

If  $A$ and $B$ are commutative algebras where $A$ is $B$-projective, 
 then the bialgebroid $T$ is just 
the tensor algebra $A \o_B A$, $R = A$,
with Sweedler $A$-coring $\cop(x \o y) =
x \o 1 \o y$ and $\eps = \mu$.  
This particular bialgebroid has an antipode $\tau: T \to T$ given by $\tau(x \o y) = y \o x$
(cf.\ \cite{Lu, PX, KS}).  

P.\ Xu \cite{PX} defines bialgebroid using an anchor map $T \to \End R$
instead of the counit $\eps: T \to R$. 
The anchor map is a right $T$-module algebra structure on $R$ given
by 
\begin{equation}
\label{eq: ract1}
r \ract t = t^1 r t^2,
\end{equation} 
for $r \in R, t \in T$.  We will
study this and an extended right $T$-module algebra structure on $\End {}_BA$ in the next section.  
The counit  is the anchor map evaluated
at $1_R$, which is indeed the case above.

\begin{remark}
\begin{rm}
If $I$ is a finite set, a D2 extension
$A \| B$ has a left bialgebroid structure on $S =
\End {}_BA_B$ such that $A$ is left
$S$-module algebra, the left
or right endomorphism algebras are smash
products of $A$ with $S$
and $T$ is the $R$-dual bialgebroid of $S$ \cite{KS}.
In the proofs of these facts, most of the formulas in \cite{KS}
do not make sense if $I$ is an infinite
set.
\end{rm}
 \end{remark}
\section{A right $T$-module endomorphism algebra}

We continue in this section with a right  depth two extension $A \| B$
and  our notation for $T = (A \o_B A)^B$, $R= C_A(B)$, left and right D2 quasibases $t_i, u_i \in T$, $\beta_i, \gamma_i \in S$
where $i \in I$, respectively, in a index set $I$ of possibly infinite cardinality.
Given any right $R$-bialgebroid $T$, recall that a right $T$-module algebra
$A$ is an algebra in the tensor category of right  $T$-modules \cite{BW, KS}.  
  
Suppose ${}_AM$ is a left $A$-module.  Let $\ME$ denote
its endomorphism ring as a module restricted to a  $B$-module: $\ME = \End {}_BM$.
There is a right action of $T$ on $\mathcal{E}$ 
given by $f \ract t = t^1 f(t^2 -)$ for $f \in \mathcal{E}$. This is
a measuring action and $\ME$ is a right $T$-module
algebra (as defined in \cite{KS, BW}), since
$$ (f \ract t\1)\circ (g \ract t\2) = \sum_i t^1 f(\gamma_i(t^2) u^1_ig(u^2_i -)) = (f\circ g) \ract t, $$
and $1_{\ME} \in \ME^T$ is a $T$-invariant, since $\id_M \ract t = \id_M \ract s_R(\eps(t))$.  
The subring of invariants $\ME^T$ in $\mathcal{E}$ is $\End {}_AM$ since $\End {}_AM \subseteq \ME^T$ is obvious,
and $\phi \in \ME^T$ satisfies for $m \in M, a \in A$:
$$ \phi(am) = \sum_i \gamma_i(a) (\phi \ract u_i)(m)
= \sum_i \gamma_i(a) \eps_T(u_i) \phi(m) = a  \phi(m).$$
With similar arguments for left D2 quasibase, we have established:

\begin{theorem}
If $B \to A$ is right or left D2 and ${}_AM$ is module, then $\End {}_BM$
is a right $T$-module algebra with invariant subalgebra
$\End {}_AM$.
\end{theorem}

By specializing $M= A$, we obtain 

\begin{cor}
If $A \| B$ is D2, then $\End {}_BA$ is a right $T$-module algebra
with  invariant subalgebra $\rho(A)$ and right $T$-module subalgebra $\lambda(R)$.  
\end{cor}  
\begin{proof}
For any $r \in R$, we have $\lambda(r) \ract t = \lambda (r \ract t)$ wrt.\ the right action of $T$ on $R$ in eq.~(\ref{eq: ract1}) in the previous section. Of course, $\Hom ({}_AA, {}_AA) \cong \rho(A)$ where we fix the notation for right multiplication, $\rho(a)(x) = xa$ (all $a,x \in A$). 
\end{proof}

The right $T$-module $\End {}_BA$
is identifiable with composition of endomorphism and homomorphism under the ring isomorphism $T \cong \End {}_AA \o_B A_A$ and the
$A$-$A$-bimodule isomorphism 
$\End {}_BA \cong \Hom ({}_AA \o_B A, {}_AA)$ via
$f \mapsto (x \o y \mapsto xf(y))$.  We leave this remark
as an exercise.


\section{Main theorem characterizing Galois extension}
Given any right $R$-bialgebroid $T$, recall that a right $T$-comodule algebra
$A$ is an algebra in the category of right $T$-comodules
\cite{BW}. If $B$ denotes its subalgebra of coinvariants $A^{\rm co \, T}$, which are
the elements $\delta: x \mapsto x \o_R 1_T$ under the coaction, we
say $A \| B$ is right $T$-Galois if the canonical mapping
$\beta: A \o_B A \to A \o_R T$ given by $\beta(x \o y) = x y\0 \o y\1$
is bijective.  Note that any $r \in R$ and $b \in B$ necessarily commute
in $A$, since
the coaction is monic and 
$$\delta (rb) = \delta(r)\delta(b) =  b \o_R s_R(r) = \delta (br) .$$

Among other things, we  show in the theorem that if $A \| B$ is right
depth two, then $A$ is a right $T$-comodule algebra
and the isomorphism $A \o_B A \cong A \o_R T$ projects to
the Galois mapping via $A \o_B A \to A \o_{A^{\rm co \, T}} A$.
If moreover the natural module $A_B$ is faithfully flat (apply to
eq.~\ref{eq: ff} below) or balanced, 
i.e., the map $\rho$: $B \to \End {}_EA$ is surjective
 where $E = \End A_B$, then $B = A^{\rm co \, T}$.       

\begin{theorem}
\label{th-main}
An algebra extension $A \| B$ is right D2 and right balanced
if and only if $A \| B$ is a right $T$-Galois extension for
some left projective right bialgebroid $T$ over some algebra $R$.  
\end{theorem}
\begin{proof}
($\Leftarrow$)  
Since ${}_RT$ is projective, ${}_RT \oplus * \cong {}_R R^{(I)}$
for some set $I$.  Then $A \o_R T \oplus * \cong A^{(I)}$
as $A$-$B$-bimodules, since $R$ and $B$ commute in $A$
and $A \o_R R^{(I)} \cong (A \o_R R)^{(I)}$.  But the Galois isomorphism $A \o_B A \stackrel{\cong}{\to} A \o_R T$,
$\beta(x \o y) = xy\0 \o y\1$
is an $A$-$B$-bimodule isomorphism, hence
$A \| B$ is right D2.

To see that the natural map $\rho: B \to \End {}_EA$ is surjective,
we let $F \in \End {}_EA$.  Then for each $a \in A$, left multiplication
$\lambda_a \in E$, whence $F \circ \lambda_a = \lambda_a \circ F$.
It follows that $F = \rho_x$ where $x = F(1)$.  Since
$B = A^{\rm co \, T}$, it suffices 
to show  that $x\0 \o x\1 = x \o 1$ under the coaction.  For this
we pause for a lemma.

Lemma:  Let $R$ be an algebra with modules $M_R$ and ${}_RV$
where $V$ is projective with dual bases $w_i \in V$, $f_i \in \Hom ({}_RV, {}_RR) = {}^*V$ for some possibly infinite cardinality index set $i \in I$.  
If for some $m_ j \in M, v_j \in V$ and finite index set $J$, we have$\sum_{j \in J} m_j \phi(v_j) = 0$ for each $\phi \in {}^*V$,
then $\sum_{j \in J} m_j \o_R v_j = 0$.  This statement
follows of course by substitution of $\sum_{i \in I} f_i(v_j)w_i$
for each $j \in J$.  

To see that $x \o 1 - x\0 \o x\1 = 0$, we define for each $\nu \in \Hom ({}_RT, {}_RR)$, the right $B$-endomorphism $\overline{\nu} \in E$
by $\overline{\nu}(a) = a\0 \nu(a\1)$.  Since also $\rho_r \in E$ for
$r = \nu(1_T) \in R$, we compute:
$$ x \nu(1_T) = \rho_{\nu(1_T)} F(1_A) = F(\overline{\nu}(1_A)) = 
\overline{\nu}(F(1_A)) = x\0 \nu(x\1). $$
By lemma then $x \in B$, so that $A_B$ is a balanced module.  

($\Rightarrow$) If $A \| B$ is right D2, we have explicit 
formulas in the previous section for $T = (A \o_B A)^B$
as a left $R$-projective right bialgebroid over $R = C_A(B)$.  
Define a coaction $\delta: A \to A \o_R T$ by \begin{equation}
\label{eq: coaction}
\delta(a) = \sum_{i \in I}
\gamma_i(a) \o_R u_j.
\end{equation}
  We claim that $A$ is a right $T$-comodule algebra,
an argument similar to \cite[5.1]{LK2005} but with infinite index set,
and postpone a sketch of the proof for two paragraphs.

It is clear that $B \subseteq A^{\rm co \, T}$, since $\gamma_i(b) = 
b\gamma_i(1)$ and $\sum_i \gamma(1)u_i = 1_T$.  For the reverse
inclusion, suppose $x \in A$ such that $x\0 \o x\1 = x \o 1_T$.
Note the isomorphism
\begin{equation}
\label{eq: ice}
A \o_R T \stackrel{\cong}{\longrightarrow} A \o_B A 
\end{equation}
given by $a \o t \mapsto at^1 \o t^2$,
with inverse
\begin{equation}
\label{eq: beta}
\beta(x \o y) = \sum_i x\gamma_i(y) \o_R u_i, \ \ \beta: A \o_B A
\stackrel{\cong}{\longrightarrow} A \o_R T .
\end{equation}
 Since $\sum_i \gamma_i(x) \o u_i = x \o 1_T$,
the image under $\beta^{-1}$ is 
\begin{equation}
\label{eq: ff}
1_A \o_B x = x \o_B 1_A.
\end{equation}
Given any $f \in E$, apply $\mu \circ (f \circ \lambda_a \o \id_A)$
to this equation obtaining $f(a)x =f(ax) $ for each $a \in A$.
Then $\rho_x \circ f = f \circ \rho_x$, whence $\rho_x \in \End {}_EA$.
Then $\rho_x \in \rho(B)$ since $A_B$ is balanced.  Hence $x \in B$.  

The Galois condition on the algebra extension $A \| B$ follows immediately from the fact that $\beta$ in
eq.~(\ref{eq: beta}) is an isomorphism. Indeed using the isomorphism
$\beta^{-1}$ as an identification between $A \o_R T$ and $A \o_B A$
is the easiest way to show $\delta$ defines a right $T$-comodule structure
on $A$. 
  
The conditions that $A$ must meet to be a
right $T$-comodule algebra are
\begin{enumerate}
\item an algebra homomorphism $R \to A$;
\item a right $T$-comodule structure $(A, \delta)$:
$\delta$ is right $R$-linear, $a\0 \eps(a\1) = a$ for all $a \in A$,
$(\id_A \o \cop) \circ \delta = (\delta \o \id_T) \circ \delta$;
\item $\delta(1_A) = 1_A \o 1_T$;
\item for all $r \in R, a \in A$, $r a\0 \o_R a\1 = a\0 \o_R t_R(r) a\1$;
\item $\delta(xy) = x\0 y\0 \o_R x\1 y\1$ for all $x,y \in A$.
\end{enumerate}
The following is a sketch of the proof, the details being left to the reader. 
For $R \to A$ we take the inclusion $C_A(B) \into A$.  Note that
$\delta(ar) = a\0 \o_R a\1 s_R(t)$ since both expressions map into $1 \o_B ar$ under $\beta^{-1}: A \o_R T \stackrel{\cong}{\longrightarrow} A \o_B A$.  Next we note that $A$ is counital since $\sum_i \gamma_i(a)u^1_i u^2_i = a$. The coaction is coassociative on any $a \in A$ since both expressions 
map into $1_A \o 1_A \o a$ under the isomorphism 
\begin{equation}
\label{eq: isom}
A \o_R T \o_R T\stackrel{\cong}{\longrightarrow} A \o_B A \o_B A, \ \ \
a \o t \o u \longmapsto at^1 \o t^2 u^1 \o u^2. 
\end{equation}
The expressions in the last two items map bijectively via $\beta^{-1}$ into
$r \o a$ and $1 \o xy$ in $A \o_B A$, respectively, so the equalities hold. 
\end{proof}

The following by-product of the proof above is a characterization
of right (similarly left) depth two
in terms of $T$.

\begin{cor}
Let $A \| B$ be an algebra extension
with $T = (A \o_B A)^B$ and $R = C_A(B)$.
The extension $A \| B$ is right D2 if and only if $A \o_R T \cong A \o_B A$
via $a \o_R t \mapsto at^1 \o_B t^2$
and the module ${}_RT$ is projective.
\end{cor}

The main theorem is most interesting for
subalgebras with small centralizers.
An example of what can happen for large centralizers:  the theorem shows
that \textit{any} field extension $K \supseteq F$ is $T$-Galois, since the underlying vector space of the $F$-algebra $K$ is 
 free, therefore balanced, 
and any algebra over a field is depth two. 
The bialgebroid $T$ in this case is
remarked on after Theorem~\ref{th-bi}. 

The paper \cite{LK2005} sketches how
the main theorem in this paper would
 extend the main
theorem in \cite{KN} for 
extensions with trivial centralizer as follows.
We call an algebra extension $A \| B$ 
semisimple-Hopf-Galois if $H$ is
a semisimple Hopf algebra, $A$
is an $H$-comodule algebra
with coinvariants $B$ and the Galois
mapping $A \o_B A \to A \o H$ is bijective \cite{Mo}.  Recall
that an algebra extension $A \| B$
is a Frobenius extension if $A_B$
if f.g.\ projective and $A \cong
\Hom (A_B, B_B)$ as natural
$B$-$A$-bimodules.  Left and right depth
two are equivalent conditions on a Frobenius
extension \cite{KS}. Recall too
that an algebra extension $A \| B$ is separable
if the multiplication $\mu: A \o_B A \to A$ is
a split $A$-$A$-epi.   

\begin{cor}
Suppose $A \| B$ is a Frobenius extension of $k$-algebras with trivial centralizer $R = 1_A \cdot k$ and $k$ a field of characteristic zero.  
Then $A \| B$ is semisimple-Hopf-Galois if and only if $A \| B$ is a  separable and depth two extension.
\end{cor}

Also various pseudo-Galois and almost-Galois extensions over groups,
Hopf algebras or  weak Hopf algebras are depth two, balanced 
extensions, and so Galois extensions with
respect to bialgebroids.  The following
corollary is an example using Hopf algebras, although the corollary
may be stated more generally for bialgebroids by using
the proof of $\Leftarrow$ above, which stays valid if the Galois mapping $\beta: A \o_B A \to A \o_R T$ is weakened from isomorphism to split $A$-$B$-monomorphism. 

\begin{cor}
Suppose $H$ is a Hopf algebra and $A \| B$ is a right $H$-extension.  If the Galois
mapping $\beta$ is a split $A$-$B$-monomorphism, then $A \| B$ is a right $(A \o_B A)^B$-Galois extension,
where $(A \o_B A)^B$ is the bialgebroid $T$ 
over $C_A(B)$ studied in section~2.  
\end{cor}
  
This corollary fits in with the current  study of weak Hopf-Galois extensions in case 
the centralizer $C_A(B)$ is a separable algebra over a field, whence the bialgebroid $T$ is a weak bialgebra \cite[Prop.\ 7.4]{KS}.  


\end{document}